\documentclass{amsart}
\usepackage{cases}
\usepackage{graphicx}
\usepackage{amssymb}
\usepackage{mathrsfs}

\numberwithin{equation}{section}

\newcommand{\bqn}{\begin{equation}}
\newcommand{\eqn}{\end{equation}}

\begin{document}
\title[]
{Ramanujan's cubic transformation and generalized modular equation}
\author{Miao-Kun Wang}

\address{%
Department of Mathematics, Huzhou Teachers College, Huzhou,
Zhejiang, 313000, China}\

\email{wmk000@126.com}

\thanks{This work was supported by the Natural Science
Foundation of China (Grant Nos. 11071069, 11171307) and PhD Students
Innovation Foundation of Hunan Province (CX2012B153).}
\author{Yu-Ming Chu}

\address{%
Department of Mathematics, Huzhou Teachers College,
    Huzhou, Zhejiang, 313000, China }

\email{chuyuming@hutc.zj.cn}

\author{Yue-Ping Jiang}

\address{%
College of Mathematics and Econometrics, Hunan University, Changsha,
Hunan, 410082, China}

\email{ypjiang731@163.com}



\begin{abstract}
\footnotesize  We study the quotient of hypergeometric functions
\begin{equation*}
\mu_{a}^*(r)=\frac{\pi}{2\sin{(\pi
a)}}\frac{F(a,1-a;1;1-r^3)}{F(a,1-a;1;r^3)} \quad (r\in(0,1))
\end{equation*}
in the theory of Ramanujan's generalized modular equation for
$a\in(0,1/2]$, find an infinite product formula for $\mu_{1/3}^*(r)$
by use of the properties of $\mu_{a}^*(r)$ and Ramanujan's cubic
transformation. Besides, a new cubic transformation formula of
hypergeometric function is given, which complements the Ramanujan's
cubic transformation.
\end{abstract}

\maketitle

\section{Introduction}
For real numbers $a$, $b$ and $c$ with $c\neq 0, -1, -2, \cdots$,
the Gaussian hypergeometric function is defined by
\begin{equation}
F(a,b;c;x)={}_{2}F_{1}(a,b;c;x)=\sum_{n=0}^{\infty}\frac{(a,n)(b,n)}{(c,n)}\frac{x^n}{n!}
\end{equation}
for $x\in(-1,1)$, where $(a,n)$ denotes the shifted factorial
function
$$(a,n)=a(a+1)(a+2)(a+3)\cdots(a+n-1)$$
for $n=1,2,\cdots$, and $(a,0)=1$ for $a\neq 0$. In particular,
$F(a,b;c;x)$ is called zero-balanced if $c=a+b$.

It is well known that $F(a,b;c;x)$ has many important applications
in various fields of the mathematical and natural sciences, and many
classes of special function in mathematical physics are particular
cases of this function. For these, and properties of  $F(a,b;c;x)$
see [1, 2, 4-6, 9, 16, 22, 27, 31, 39]. Here we recall one of the
most important properties of $F(a,b;c;x)$, the Ramanujan's cubic
transformation,
\begin{equation}
F(\frac{1}{3},\frac{2}{3};1;1-\left(\frac{1-r}{1+2r}\right)^3)=(1+2r)F(\frac{1}{3},\frac{2}{3};1;r^3)
\end{equation}
or
\begin{equation}
F(\frac{1}{3},\frac{2}{3};1;\left(\frac{1-r}{1+2r}\right)^3)=\frac{1+2r}{3}F(\frac{1}{3},\frac{2}{3};1;1-r^3),
\end{equation}
which was raised by S. Ramanujan in his unpublished notebooks. In
1989, J. M. and  P. M. Borwein [19] provided a new proof of equation
(1.2) or (1.3).

\medskip
{\sc Theorem 1.1} (see [19]). \emph{Let $a,b>0$,
\begin{align}
&a_{n+1}:=\frac{a_{n}+2b_{n}}{3}, &&a_{0}:=a,\\
&b_{n+1}:=\sqrt[3]{\frac{b_{n}(a_{n}^2+a_{n}b_{n}+b_{n}^2)}{3}},&&
b_{0}:=b.
\end{align}
Then when $a=1$ and $b=x\in(0,1)$, the common limit $F(1,x)$ of
$\{a_{n}\}$ and $\{b_{n}\}$ is given by
\begin{equation*}
\frac{1}{F(1,x)}=\sum\limits_{n=0}^{\infty}\frac{(1/3,n)(2/3,n)}{(n!)^2}(1-x^3)^n=F(\frac{1}{3},\frac{2}{3};1;1-x^3).
\end{equation*}}

For $a\in(0,1/2]$, $r\in(0,1)$, $p>0$, Ramanujan's generalized
modular equation with signature $1/a$ and degree $p$ is defined by
\begin{equation}
\frac{F(a,1-a;1;1-s^2)}{F(a,1-a;1;s^2)}=p\frac{F(a,1-a;1;1-r^2)}{F(a,1-a;1;r^2)}.
\end{equation}
Making use of the decreasing homeomorphism $\mu_{a}:(0,1)\rightarrow
(0,\infty)$ defined by
\begin{equation}
\mu_{a}(r)=\frac{\pi}{2\sin{(\pi
a)}}\frac{F(a,1-a;1;1-r^2)}{F(a,1-a;1;r^2)},
\end{equation}
we can rewrite (1.6) as
\begin{equation}
\mu_{a}(s)=p\mu_{a}(r),\quad 0<r<1.
\end{equation}
The solution of (1.8) is given by
\begin{equation*}
s\equiv\varphi_{K}(a,r)={\mu_{a}}^{-1}(\mu_{a}(r)/K),\quad K=1/p.
\end{equation*}
In the particular case $a=1/2$, Ramanujan's generalized modular
equation (1.6) reduces to the classical case, and the modular
functions $\mu_{a}(r)$ and $\varphi_{K}(a,r)$ become $\mu(r)$ and
$\varphi_{K}(r)$, respectively.

In this note, for convenience, we denote $r^*=\sqrt[3]{1-r^3}$,
\begin{equation}
\mu_{a}^*(r)=\mu_{a}(r^{3/2})=\frac{\pi}{2\sin{(\pi
a)}}\frac{F(a,1-a;1;1-r^3)}{F(a,1-a;1;r^3)},
\end{equation}
\begin{equation}
\varphi_{K}^*(a,r)=\left[\varphi_{K}(a,r^{3/2})\right]^{2/3}={\mu_{a}^*}^{-1}(\mu_{a}^*(r)/K)
\end{equation}
and
\begin{equation}
\mu^*(r)=\mu_{1/3}^*(r), \quad
\varphi_{K}^*(r)=\varphi_{K}^*(1/3,r).
\end{equation}
Then from (1.2), (1.3) and (1.9)-(1.11) we conclude that
\begin{equation}
\mu_{a}^*(r)\mu_{a}^*(r^*)=\frac{{\pi}^2}{4\sin^2(\pi a)},
\end{equation}
\begin{equation}
\mu^*(r)=3\mu^*\left(\frac{\sqrt[3]{9r(1+r+r^2)}}{1+2r}\right),\quad
\mu^*(r)=\frac{1}{3}\mu^*\left(\frac{1-r^*}{1+2r^*}\right),
\end{equation}
\begin{equation}
\mu^*(r)\mu^*\left(\frac{1-r}{1+2r}\right)=\pi^2,
\end{equation}
\begin{equation}
\varphi_{3}^*(r)=\frac{\sqrt[3]{9r(1+r+r^2)}}{1+2r},\quad
\varphi_{1/3}^*(r)=\frac{1-r^*}{1+2r^*},\quad
\varphi_{3}^*(r)^3+\varphi_{1/3}^*(r^*)^3=1.
\end{equation}

It follows from (1.9) and (1.10) that in order to study the modular
functions $\mu_{a}(r)$ and $\varphi_{K}(a,r)$, we only need to
consider the functions $\mu_{a}^*(r)$ and $\varphi_{K}^*(a,r)$.

\medskip
As is known to all, Ramanujan's cubic transformation and generalized
modular equation have been developed for over a century. In 1900s,
S. Ramanujan studied extensively $F(a,b;c;x)$ and the modular
equation (1.6), and gave a lot of statements concerning them in his
unpublished notebooks [28-30], but no original proof have remained.
Later, Ramanujan's theories have been developed by many authors,
such as J.M. and P.B. Borwein [17-20], K. Venkatachaliengar [34] and
B. C. Berndt [11-15].

The greatest advances toward establishing Ramanujan's theories have
been made by J. M. and P. B. Borwein [20]. In searching for
analogues of the classical arithmetic-geometric mean of Gauss, J. M.
and P. B. Borwein discovered an elegant cubic analogue, namely,
(1.4) and (1.5). Thus a cubic transformation formula (1.2) or (1.3)
for $F(1/3,2/3;1;x)$ was derived. In fact, equations (1.2) and (1.3)
can be found on page 258 of Ramanujan's second notebook [28], and
they were rediscovered by the Borweins.

In 1995, B. C. Berndt, S. Bhargava, and F. G. Garvan published a
landmark paper [15] in which they studied the generalized modular
equation (1.6) with $p$ an integer. For several rational values of
$a$ such as $a=1/3, 1/4, 1/6$ and prime $p$ (e.g. $p=2, 3, 5, 7,
\cdots$), they were able to give proofs for numerous algebraic
identities stated by S. Ramanujan in his unpublished notebooks.
Meanwhile, a generalization of Ramanujan's cubic transformation for
$F(1/3,2/3;1;x)$ was given. After the publication of [15] many
papers have been written on modular equations [3, 10, 26, 32].

A new connection between geometric function theory and Ramanujan's
theory was derived by M. Vuorinen in [35]. He found that the
functions $\varphi_{K}(r)$ and $\mu(r)$, as the Hersch-Pfluger
distortion function and the plane Gr\"{o}tzsch ring function, play
an important role in the theory of quasiconformal maps. Then the
functions $\varphi_{K}(r)$ and $\mu(r)$, and their generalizations
$\varphi_{K}(a,r)$ and $\mu_{a}(r)$ have been the subject of
intensive research. In particular, many remarkable inequalities for
them can be found in the literature [3, 7, 21, 24-26, 36, 38].
Especially, G. D. Anderson, S.-L. Qiu, M. K. Vamanamurithy and M.
Vuorinen [3] established several analytic properties for
$\varphi_{K}(a,r)$, as applications, some estimates to the solution
of generalized modular equation (1.6) were obtained. Recently, G.-D.
Wang, X.-H. Zhang and Y.-M. Chu [36] found the relation between the
modular function $\varphi_{K}(a,r)$ and $\mu_{a}(r)$, and proved
that, for  $a(r)$ is a real function defined on $(0,1)$,
$r\in(0,1)$ and $K\in(1,\infty)$, inequality
\begin{equation*}
\varphi_{1/K}(a,r)>r^K\exp\left\{(1-K)a(r)\right\}
\end{equation*}
holds if and only if $a(r)\geq \mu_{a}(r)+\log{r}$. Equivalently, by
(1.9) and (1.10),
\begin{equation}
\varphi_{1/K}^*(a,r)>r^K\exp\left\{\frac{2}{3}(1-K)a^*(r)\right\}
\end{equation}
if and only if $a^*(r)\geq \mu_{a}^*(r)+(3/2)\log{r}$, where
$a^*(r)$ is also a real functions defined on $(0,1)$.

\medskip
The main purpose of this paper is to find a infinite-product
representation for $\mu^*(r)$(or $\mu_{1/3}(r)$) which only contains
$r$, and to extend representation to the function $\mu_{a}^*(r)$. We
shall prove the following Theorem 1.2.

\medskip {\sc Theorem 1.2.} \emph{For $a\in(0,1/2]$ and $r\in(0,1)$, let
$r_{0}=r^*=\sqrt[3]{1-r^3}$,
$r_{1}=\varphi_{3}^*(r^*)=\sqrt[3]{9r^*(1+r^*+{r^*}^2)}/(1+2r^*),
\cdots$, and
\begin{equation}
r_{n}=\varphi_{3}^*(r_{n-1})=\frac{\sqrt[3]{9r_{n-1}(1+r_{n-1}+{r_{n-1}}^2)}}{1+2r_{n-1}}=\varphi_{3^n}^*(r^*).
\end{equation}
Then
\begin{align}
&\prod\limits_{n=0}^{\infty}\left[(1+2r_{n})(1+r_{n}+r_{n}^2)\right]^{\frac{3^{-n}}{2}}\leq
\exp(\mu_{a}^*(r)+\frac{3}{2}\log{r})\\
&\leq
\frac{1}{\sqrt{27}}\exp\left(R(a)/2\right)\prod\limits_{n=0}^{\infty}\left[(1+2r_{n})(1+r_{n}+r_{n}^2)\right]^{\frac{3^{-n}}{2}}\nonumber
\end{align}
for $a\in(0,1/3]$, and the revered inequality
\begin{align}
&\frac{1}{\sqrt{27}}\exp\left(R(a)/2\right)\prod\limits_{n=0}^{\infty}\left[(1+2r_{n})(1+r_{n}+r_{n}^2)\right]^{\frac{3^{-n}}{2}}\leq
\exp(\mu_{a}^*(r)+\frac{3}{2}\log{r})\\
&\leq
\prod\limits_{n=0}^{\infty}\left[(1+2r_{n})(1+r_{n}+r_{n}^2)\right]^{\frac{3^{-n}}{2}}\nonumber
\end{align}
holds for $a\in[1/3,1/2]$. Moreover, each equality in (1.18) and
(1.19) is reached if and only if $a=1/3$. In particular, for all
$r\in(0,1)$,
\begin{equation}
\exp(\mu^*(r)+\frac{3}{2}\log{r})=\exp(\mu_{1/3}(r^{3/2})+\frac{3}{2}\log{r})=\prod\limits_{n=0}^{\infty}\left[(1+2r_{n})(1+r_{n}+r_{n}^2)\right]^{\frac{3^{-n}}{2}}
\end{equation}
or
\begin{equation}
\mu^*(r)+\frac{3}{2}\log{r}=\mu_{1/3}(r^{3/2})+\frac{3}{2}\log{r}=\frac{1}{2}\sum\limits_{n=0}^{\infty}\frac{1}{3^{n}}\log\left[(1+2r_{n})(1+r_{n}+r_{n}^2)\right].
\end{equation}}

Theorem 1.2 and inequality (1.16) lead to the following corollary.

\medskip
{\sc Corollary 1.3.} \emph{For $a\in(0,1/2]$, $r\in(0,1)$ and
$K\in(1,\infty)$, let $r_{n}$ defined as in Theorem 1.2, and
$\Sigma\equiv\sum_{n=0}^{\infty}{3^{-n}}\log\left[(1+2r_{n})(1+r_{n}+r_{n}^2)\right]$.
Then
\begin{equation*}
\varphi_{1/K}^*(a,r)>r^{K}\exp\left\{\frac{(1-K)}{3}\left(R(a)-\log{27}
+\Sigma\right)\right\}
\end{equation*}
for $a\in(0,1/3]$ and $r\in(0,1)$, and
\begin{equation*}
\varphi_{1/K}^*(a,r)>r^{K}\exp\left\{\frac{\Sigma(1-K)}{3}\right\}.
\end{equation*}
for $a\in[1/3,1/2)$ and $r\in(0,1)$.}

\medskip
Another purpose of this paper is to complement Theorem 1.1. From
(1.4) and (1.5) we clearly see that the iteration is positively
homogeneous so is the limit function $F(a,b)$. Without loss of
generality, we only consider two cases: (A) $a=1,b=x\in(0,1)$; (B)
$a=x\in(0,1),b=1$. Theorem 1.1 gives the limit function of case A,
while the following Theorem 1.4 presents the answer of case B.

\medskip
{\sc Theorem 1.4.} \emph{Let $a,b>0$,
\begin{align*}
&a_{n+1}:=\frac{a_{n}+2b_{n}}{3}, &&a_{0}:=a,\\
&b_{n+1}:=\sqrt[3]{\frac{b_{n}(a_{n}^2+a_{n}b_{n}+b_{n}^2)}{3}},&&
b_{0}:=b.
\end{align*}
Then when $a=x\in(0,1)$ and $b=1$, the common limit $F(x,1)$ of
$\{a_{n}\}$ and $\{b_{n}\}$ is
\begin{equation*}
\frac{1}{F(x,1)}=\sum\limits_{n=0}^{\infty}\frac{(1/3,n)^2}{(n!)^2}(1-x^3)^n=F(\frac{1}{3},\frac{1}{3};1;1-x^3).
\end{equation*}
In particular, for $x\in(0,1)$, then
\begin{equation}
F\left(\frac{1}{3},\frac{1}{3};1;1-x^3\right)=\sqrt[3]{\frac{3}{x^2+x+1}}F\left(\frac{1}{3},\frac{1}{3};1;\frac{(1-x)^3}{9(x^2+x+1)}\right).
\end{equation}}

\medskip
The methods of the proofs of Theorems 1.2 and 1.4  primarily come
from S.-L. Qiu and M. Vuorinen in [24], and J. M. and P. B. Borwein
in [19], respectively.

\bigskip
\section{Preliminary results}
\setcounter{equation}{0}

In this section, we study some monotonicity properties of the
modular function $\mu_{a}^*(r)$, which will be used in the proof of
Theorem 1.2. But first, we recall some known results for the
function $F(a,b;c;x)$.

It is well known that the properties of the hypergeometric functions
are closely related to those of the gamma function $\Gamma(x)$, the
psi function $\Psi(x)$, and the beta function $B(x,y)$. For positive
numbers $x$ and $y$, these functions are defined by
\begin{equation}
\Gamma(x)=\int_{0}^{\infty}e^{-t}t^{x-1}dt,\quad
\Psi(x)=\frac{\Gamma'(x)}{\Gamma(x)}, \quad
B(x,y)=\frac{\Gamma(x)\Gamma(y)}{\Gamma(x+y)},
\end{equation}
respectively (cf. [39]). It is well known that the gamma function
satisfies the difference equation
\begin{equation*}
\Gamma(x+1)=x\Gamma(x)
\end{equation*}
if $x$ is nonpositive integer and has the so-called reflection
property
\begin{equation*}
\Gamma(x)\Gamma(1-x)=\frac{\pi}{\sin{\pi x}}=B(x,1-x)
\end{equation*}
if $x$ is not an integer. We shall also need the function
\begin{equation}
R(a,b)=-2\gamma-\Psi(a)-\Psi(b), \quad R(a)=R(a,1-a), \quad
R(1/3,2/3)=\log{27},
\end{equation}
where $\gamma$ is the Euler-Mascheroni constant defined by
\begin{equation*}
\gamma=\lim_{n\rightarrow
\infty}\left(\sum\limits_{k=1}^{n}\frac{1}{k}-\log{n}\right)=0.577215\cdots.
\end{equation*}
By [25, Lemma 2.14(2)], $R(a)$ is strictly decreasing in
$a\in(0,1/2]$. Thus $R(a)>\log27$ for $a\in(0,1/3)$, and
$R(a)<\log27$ for $a\in(1/3,1/2]$.

One important tool we shall need in our work is the following
Ramanujan's derivative formula [12, Corollary, p.86]
\begin{equation}
\frac{d}{dx}\left[\frac{F(a,1-a;1;1-x)}{F(a,1-a;1;x)}\right]=-\frac{\sin(\pi
a)}{\pi x(1-x)F(a,1-a;1;x)^2}
\end{equation}
for $a,x\in(0,1)$. Then from (2.3) we immediately get the derivative
of $\mu_{a}^*(r)$ with respect to $r$: for $a \in(0,1/2)$ and $r
\in(0, 1)$,
\begin{equation}
\frac{d\mu_{a}^*(r)}{dr}=-\frac{3}{2}\cdot\frac{1}{r(1-r^3)F(a,1-a;1;r^3)^2}.
\end{equation}

Another important tool in our work is the following Ramanujan's
cubic transformation inequalities for zero-balanced hypergeometric
function.

\medskip
{\sc Theorem 2.1} (see [37, Theorem 2.4]). \emph{Let $B(a,b)$ and
$R(a,b)$ are defined as in (2.1) and (2.2), respectively. Then for
$(a,b)\in \{(a,b)|a,b>0,a+b\leq 1,ab-2(a+b)/9\leq 0\}$, inequality
\begin{align}
0\leq&
(1+2r)F(a,b;a+b;r^3)-F(a,b;a+b;\frac{9r(1+r+r^2)}{(1+2r)^3})\\
\leq& \frac{2(R(a,b)-\log{27})}{B(a,b)}\nonumber
\end{align}
holds for all $r\in(0,1)$. Also, for $(a,b)\in \{(a,b)|a,b>0,a+b\geq
1,ab-2(a+b)/9\geq 0\}$,
\begin{align}
0\leq&
F(a,b;a+b;\frac{9r(1+r+r^2)}{(1+2r)^3})-(1+2r)F(a,b;a+b;r^3)\\
\leq& \frac{2(\log{27}-R(a,b))}{B(a,b)}.\nonumber
\end{align}}

Other important tool in the rest of this paper is the following
Lemma 2.2.

\medskip
{\sc Lemma 2.2} (see [23, Lemma 2.1]). \emph{Suppose that the power
series $f(x)=\sum_{n=0}^{\infty}a_{n}x^{n}$ and
$g(x)=\sum_{n=0}^{\infty}b_{n}x^{n}$ have the radius of convergence
$r>0$, $a_{n}\in \mathbb{R}$ and $b_{n}>0$ for all
$n\in\{0,1,2,\cdots\}$. Let $h(x)={f(x)}/{g(x)}$, then the sequence
$\{a_{n}/b_{n}\}_{n=0}^{\infty}$ is (strictly) increasing
(decreasing), then $h(x)$ is also (strictly) increasing (decreasing)
on $(0,r)$.}

\medskip
Motivated by S. Simi\'{c} and M. Vuorinen [33], and \'{A}. Baricz
[8], we will employ  Lemma 2.2 to present some Ramanujan's cubic
transformation inequalities for Gaussian hypergeometric functions,
Kummer hypergeometric functions, generalized Bessel functions and
for general power series (Theorems 2.3 and 2.4, Corollary 2.5).
These results complement some results in [37], and also will be used
in the proof of Theorem 2.7.

\medskip
{\sc Theorem 2.3.} \emph{Let $a,b,c\in\mathbb{R}$ such that $c$ is
not a negative integer or zero and consider the function
$Q:(0,1)\rightarrow (0,\infty)$, defined by
$$Q(x)=F(a,b;c;x)/F(1/3,2/3;1;x).$$
Then the following assertions are true:
\begin{enumerate}
  \item If $a+b\geq c$, $9ab/2\geq \max\{1,c\}$, then $Q(x)$ is
  increasing, and consequently
\begin{equation}
F(a,b;c;\frac{9r(1+r+r^2)}{(1+2r)^3})\geq (1+2r)F(a,b;c;r^3),
\end{equation}
\begin{equation}
F(a,b;c;\left(\frac{1-r}{1+2r}\right)^3)\leq
\frac{1+2r}{3}F(a,b;c;1-r^3)
\end{equation}
hold for each $r\in(0,1)$.
  \item If $a+b\leq c$, $9ab/2\leq \min\{1,c\}$, then $Q(x)$ is
  decreasing, and consequently
\begin{equation}
F(a,b;c;\frac{9r(1+r+r^2)}{(1+2r)^3})\leq (1+2r)F(a,b;c;r^3),
\end{equation}
\begin{equation}
F(a,b;c;\left(\frac{1-r}{1+2r}\right)^3)\geq
\frac{1+2r}{3}F(a,b;c;1-r^3)
\end{equation}
hold for each $r\in(0,1)$.
\end{enumerate}}

\medskip
{\emph Proof.} Since $Q(x)$ can be written as
\begin{equation*}
Q(x)=\frac{F(a,b;c;x)}{F(\frac{1}{3},\frac{2}{3};1;x)}=\frac{\sum_{n=0}^{\infty}\frac{(a,n)(b,n)}{(c,n)}\cdot\frac{x^n}{n!}}
{\sum_{n=0}^{\infty}\frac{(1/3,n)(2/3,n)}{(1,n)}\cdot\frac{x^n}{n!}},
\end{equation*}
by Lemma 2.2, we know that the monotonicity of $Q$ depends on the
monotonicity of the sequence $\{\alpha_{n}\}$, defined by
\begin{equation*}
\{\alpha_{n}\}=\frac{(a,n)(b,n)}{(c,n)}\cdot\frac{(1,n)}{(1/3,n)(2/3,n)}.
\end{equation*}
Note that
\begin{equation*}
\frac{\alpha_{n+1}}{\alpha_{n}}=\frac{(n+a)(n+b)(n+1)}{(n+c)(n+1/3)(n+2/3)}\geq
1
\end{equation*}
if and only if
\begin{equation*}
A_{n}=(a+b-c)n^2+(a+b-c+ab-\frac{2}{9})n+ab-\frac{2}{9}\geq 0.
\end{equation*}
Thus, if $a+b\geq c$ and $9ab/2\geq \max\{1,c\}$, then $A_{n}\geq 0$
for all $n\in \{0,1,\cdots\}$, that is, the sequence
$\{\alpha_{n}\}$ is increasing, and consequently by Lemma 2.2 the
function $Q$ is increasing. Now, putting $x=x(r)=r^3$ and
$y=y(r)=9r(1+r+r^2)/(1+2r)^3$, then $0<x<y<1$ and
\begin{equation*}
\frac{F(a,b;c;r^3)}{F(\frac{1}{3},\frac{2}{3};1;r^3)}\leq
\frac{F(a,b;c;\frac{9r(1+r+r^2)}{(1+2r)^3})}{F(\frac{1}{3},\frac{2}{3};1;\frac{9r(1+r+r^2)}{(1+2r)^3})},
\end{equation*}
that is,
\begin{equation*}
F(a,b;c;r^3)\leq F(a,b;c;\frac{9r(1+r+r^2)}{(1+2r)^3})
\frac{F(\frac{1}{3},\frac{2}{3};1;r^3)}{F(\frac{1}{3},\frac{2}{3};1;\frac{9r(1+r+r^2)}{(1+2r)^3})},
\end{equation*}
which in view of (1.2) is equivalent to (2.7). Similarly, by
choosing $x=x(r)=[(1-r)/(1+2r)]^3$ and $y=y(r)=1-r^3$ we get the
inequality
\begin{equation*}
\frac{F(a,b;c;\left(\frac{1-r}{1+2r}\right)^3)}{F(\frac{1}{3},\frac{2}{3};1;\left(\frac{1-r}{1+2r}\right)^3)}\leq
\frac{F(a,b;c;1-r^3)}{F(\frac{1}{3},\frac{2}{3};1;1-r^3)},
\end{equation*}
that is, by (1.3),
\begin{equation*}
F(a,b;c;\left(\frac{1-r}{1+2r}\right)^3)\leq
\frac{1+2r}{3}F(a,b;c;1-r^3).
\end{equation*}
This proves the part (1). The proof of part (2) is similar, and thus
we omit further details. $\Box$

\medskip {\sc Remark.} If we change $r$ to $(1-r)/(1+2r)$ in
(2.7) and (2.9), then we have (2.8) and (2.10), respectively. Thus
the Ramanujan's cubic transformation inequalities (2.7) and (2.8)
are equivalent as well as the inequalities (2.9) and (2.10).

\medskip
Now, let us consider the sequence $\{\omega_{n}\}$, defined by
\begin{equation*}
\omega_{n}=\frac{(n!)^2}{(1/3,n)(2/3,n)}.
\end{equation*}
Then making use of Lemma 2.2 together with the similar argument in
Theorem 2.3 we will get a more general result of Theorem 2.3 as
follows.

\medskip
{\sc Theorem 2.4.} \emph{Suppose that the power series
$f(x)=\sum_{n=0}^{\infty}a_{n}x^n$ is convergent in all $x\in
(0,1)$, where $a_{n}\in \mathbb{R}$ for all $n\in\{0,1,\cdots\}$,
and assume that the sequence $\{a_{n}\cdot\omega_{n}\}$ is
increasing. Then the function $x\rightarrow f(x)/F(1/3,2/3;1;x)$ is
increasing on $(0,1)$, and by use of the notation
$\chi_{f}(x)=f(x^3)$ we have the Ramanujan's cubic transformation
inequalities for all $r\in(0,1)$,
\begin{equation}
\chi_{f}\left(\frac{\sqrt[3]{9r(1+r+r^2)}}{1+2r}\right)\geq
(1+2r)\chi_{f}(r).
\end{equation}
Moreover, if the sequence $\{a_{n}\cdot\omega_{n}\}$ is decreasing,
then $x\rightarrow f(x)/F(1/3,2/3;1;x)$ is decreasing on $(0,1)$,
and consequently (2.11) is reversed.}

\medskip
Applying Theorem 2.4 to the generalized Bessel function
$u_{v}:(0,\infty)\rightarrow \mathbb{R}$ and the Kummer
hypergeometric function $\Phi(p,q;\cdot): (0,\infty)\rightarrow
\mathbb{R}$, defined by
\begin{equation*}
u_{v}(x)=\sum\limits_{n=0}^{\infty}\frac{(-\frac{c}{4})^n}{(\kappa,n)}\cdot
\frac{x^n}{n!},\quad
\Phi(p,q;x)=\sum\limits_{n=0}^{\infty}\frac{(p,n)}{(q,n)}\cdot\frac{x^n}{n!},
\end{equation*}
where $v,b,c,p,q\in \mathbb{R}$, $\kappa=v+\frac{b+1}{2}\notin
\{0,-1,\cdots\}$ and $q\notin \{0,-1,\cdots\}$. Then we have the
following Corollary 2.5.

\medskip
{\sc Corollary 2.5.} \emph{Let $v,b,c,p,q\in \mathbb{R}$ such that
$\kappa\geq \max\{-1,-9c/8,-2/9-c/4\}$ and $q\geq
\max\{0,9p/2,p+7/9\}$. Then $x\rightarrow u_{v}(x)/F(1/3,2/3;1;x)$
and $x\rightarrow \Phi(p,q;x)/F(1/3,2/3;1;x)$ are decreasing on
$(0,1)$ and consequently for all $r\in(0,1)$ we get
\begin{equation*}
\chi_{u_{v}}\left(\frac{\sqrt[3]{9r(1+r+r^2)}}{1+2r}\right)\leq
(1+2r)\chi_{u_{v}}(r),
\chi_{\Phi}\left(\frac{\sqrt[3]{9r(1+r+r^2)}}{1+2r}\right)\leq
(1+2r)\chi_{\Phi}(r).
\end{equation*}}

\medskip
Next we prove the monotonicity properties and inequalities for the
modular function $\mu_{a}^*(r)$ defined as in (1.9).

\medskip
{\sc Lemma 2.6.} \emph{(1) The function $\mu_{a}^*(r)+(3/2)\log{r}$
is strictly decreasing from $(0,1)$ onto $(0,R(a)/2)$;}

\emph{(2) If $a\in(0,1/2]$, then the inequality
\begin{equation}
-\frac{3\log{3}}{2}<\mu_{a}^*(r)+\frac{1}{2}\log\left(\frac{1-r^*}{1+2r^*}\right)<\frac{R(a)}{2}
\end{equation}
holds for all $r\in(0,1)$.}

\medskip
{\emph Proof.} Part (1) directly follows from [26, Theorem 5.5(2)].

Part (2) follows from part (1) and
\begin{equation*}
\mu_{a}^*(r)+\frac{1}{2}\log\left(\frac{1-r^*}{1+2r^*}\right)=\mu_{a}^*(r)+\frac{3}{2}\log{r}-\frac{1}{2}\log(1+2r^*)(1+r^*+{r^*}^2).
\quad \Box
\end{equation*}

\medskip
{\sc Theorem 2.7.} \emph{For $a\in(0,1/2]$, let $R(a)$ be defined as
in (2.2), and $C_{1}\equiv C_{1}(a)=\min\{C, 3\}$ with
\begin{equation*}
C\equiv C(a)=\left[1+\frac{2\sin{(\pi
a)}}{\pi}(R(a)-\log{27})\right]^2.
\end{equation*}
Then we have the following statements:}

\emph{(1) For $a\in(0,1/2]$, define the function $g$ on $(0,1)$ by
\begin{equation*}
g(r)\equiv
3\mu_{a}^*\left(\frac{\sqrt[3]{9r(1+r+r^2)}}{1+2r}\right)-\mu_{a}^*(r).
\end{equation*}
Then $g$ is strictly decreasing from $(0,1)$ onto
$(0,R(a)-\log{27})$ if $a\in(0,1/3)$, is strictly increasing from
$(0,1)$ onto $(R(a)-\log{27},0)$ if $a\in(1/3,1/2]$, and
$f(r)\equiv0$ if $a=1/3$. Moreover, for $a\in(0,1/3)$ and all
$r\in(0,1)$,
\begin{equation}
\mu_{a}^*(r)<3\mu_{a}^*\left(\frac{\sqrt[3]{9r(1+r+r^2)}}{1+2r}\right)<\min\{\mu_{a}^*(r)+R(a)-\log{27},
C_{1}\mu_{a}^*(r)\}.
\end{equation}
And for $a\in(1/3,1/2)$ and $r\in(0,1)$,
\begin{equation}
\max\{\mu_{a}^*(r)+R(a)-\log{27},
C\mu_{a}^*(r)\}<3\mu_{a}^*\left(\frac{\sqrt[3]{9r(1+r+r^2)}}{1+2r}\right)<\mu_{a}^*(r).
\end{equation}}

\emph{(2) For $a\in(0,1/2]$, then the function
\begin{equation*}
f(r)\equiv \mu_{a}^*\left(\frac{1-r}{1+2r}\right)-3\mu_{a}^*(r^*)
\end{equation*}
is strictly decreasing from $(0,1)$ onto $(\log{27}-R(a),0)$ if
$a\in(0,1/3)$, and is strictly increasing from $(0,1)$ onto $(0,
\log{27}-R(a))$ if $a\in(1/3,1/2)$. Moreover, for all $a\in(0,1/3]$
and $r\in(0,1)$,
\begin{align}
&\max\left\{\frac{3\pi^2}{4\sin^2(\pi
a)}-(R(a)-\log{27})\mu_{a}^*(r),\frac{1}{C_{1}}\frac{3\pi^2}{4\sin^2(\pi
a)}\right\}\leq\\
&\mu_{a}^*(r)\mu_{a}^*\left(\frac{1-r}{1+2r}\right)\leq
\frac{3\pi^2}{4\sin^2(\pi a)}.\nonumber
\end{align}
Also, for $a\in[1/3,1/2]$ and $r\in(0,1)$,
\begin{align}
&\frac{3\pi^2}{4\sin^2(\pi a)}\leq
\mu_{a}^*(r)\mu_{a}^*\left(\frac{1-r}{1+2r}\right)\leq\\
&\min\left\{\frac{3\pi^2}{4\sin^2(\pi
a)}-(R(a)-\log{27})\mu_{a}^*(r),\frac{1}{C}\frac{3\pi^2}{4\sin^2(\pi
a)}\right\}.\nonumber
\end{align}
Equality is reached in each inequality of (2.15) and (2.16) if and
only if $a=1/3$.}

\medskip
{\emph Proof.} For part (1), if $a=1/3$, then $f(r)=0$ by (1.13).
And let $x=\sqrt[3]{9r(1+r+r^2)}/(1+2r)$, then $x^*=(1-r)/(1+2r)$
and
\begin{equation}
\frac{dx}{dr}=\frac{{x^*}^2}{3x^2}(1+2x^*)^2.
\end{equation}
It follows from Lemma 2.6(1) that
\begin{align}
&\lim\limits_{r\rightarrow 0^+}g(r)\\
=&\lim\limits_{r\rightarrow\nonumber
0^+}\left(3\mu_{a}^*(x)+\frac{9}{2}\log{x}-[\mu_{a}^*(r)+\frac{3}{2}\log{r}]+\frac{3}{2}\log{r}-\frac{9}{2}\log{x}\right)\\
=&\frac{3R(a)}{2}-\frac{R(a)}{2}-\log{27}\nonumber\\
=&R(a)-\log{27}.\nonumber
\end{align}
Clearly $\mu_{a}^*(1^-)=0$ so that
\begin{equation}
\lim\limits_{r\rightarrow 1^-}g(r)=0.
\end{equation}

Next, by differentiation and (2.4), we get
\begin{align}
g'(r)=&-\frac{9}{2}\frac{1}{x(1-x^3)F(a,1-a;1;x^3)^2}\cdot\frac{{x^*}^2}{3x^2}(1+2x^*)^2\\
&+\frac{3}{2}\frac{1}{r(1-r^3)F(a,1-a;1;r^3)^2}\nonumber\\
=&\frac{3}{2}\frac{1}{r(1-r^3)F(a,1-a;1;r^3)^2F(a,1-a;1;x^3)^2}\nonumber\\
&\times\left[F(a,1-a;1;x^3)^2-(1+2r)^2F(a,1-a;1;r^3)^2\right].\nonumber
\end{align}

Therefore, the monotonicity and range of $g$ follows from
(2.18)-(2.20) and Theorem 2.3. The first inequality and the first
upper bound in (2.13), and the first lower bound and the second
inequality in (2.14) are clear. For other inequalities in (2.13) and
(2.14), by (2.5) and (2.6) we have,  when $a\in(0,1/3]$
($a\in[1/3,1/2]$ resp.),
\begin{align}
&3\mu_{a}^*\left(\frac{\sqrt[3]{9r(1+r+r^2)}}{1+2r}\right)/\mu_{a}^*(r)\\
=&3\frac{F(a,1-a;1;r^3)}{F(a,1-a;1;1-r^3)}
\frac{F(a,1-a;1;\left(\frac{1-r}{1+2r}\right)^3)}{F(a,1-a;1;\frac{{9r(1+r+r^2)}}{(1+2r)^3})}\nonumber\\
\leq&(\geq \mbox{resp.})
\frac{(1+2r)F(a,1-a;1;r^3)}{F(a,1-a;1;1-r^3)}
\frac{F(a,1-a;1;1-r^3)+2\sin(\pi
a)(R(a)-\log{27})/\pi}{F(a,1-a;1;\frac{{9r(1+r+r^2)}}{(1+2r)^3})}\nonumber\\
\leq&(\geq \mbox{resp.}) \left[1+\frac{2\sin(\pi
a)(R(a)-\log{27})}{\pi
F(a,1-a;1;1-r^3)}\right]\left[1+\frac{2\sin(\pi
a)(R(a)-\log{27})}{\pi
F(a,1-a;1;\frac{{9r(1+r+r^2)}}{(1+2r)^3})}\right]\nonumber\\
\leq&(\geq \mbox{resp.})\left[1+\frac{2\sin(\pi
a)(R(a)-\log{27})}{\pi }\right]^2.\nonumber
\end{align}
Equality holds in each of above inequalities if and only if $a=1/3$.
on the other hand, since $x>r$, it follows from the monotonicity of
$\mu_{a}^*(r)$ with respect to $r$ on $(0,1)$  that
$\mu_{a}^*(x)<\mu_{a}^*(r)$. Hence, the remaining bounds in (2.13)
and (2.14) follow.

\medskip
For part (2), let $t=(1-r)/(1+2r)$. Then
$t^*=\sqrt[3]{9r(1+r+r^2)}/(1+2r)$ and $f(r)=-g(t)$. Hence the
assertion about $f$ follows from part (1).

It follows from (1.12), (1.15) and (2.21) that
\begin{equation}
\mu_{a}^*(r)\mu_{a}^*\left(\frac{1-r}{1+2r}\right)=\frac{3\pi^2}{4{\sin^2(\pi
a)}}\cdot
\frac{\mu_{a}^*(t)}{3\mu_{a}^*\left(\frac{\sqrt[3]{9t(1+t+t^2)}}{1+2t}\right)}\geq
\frac{1}{C_{1}}\left[\frac{3\pi^2}{4{\sin^2(\pi a)}}\right]
\end{equation}
for all $a\in(0,1/3]$ and $r\in(0,1)$, and inequality
\begin{equation}
\mu_{a}^*(r)\mu_{a}^*\left(\frac{1-r}{1+2r}\right)=\frac{3\pi^2}{4{\sin^2(\pi
a)}}\cdot
\frac{\mu_{a}^*(t)}{3\mu_{a}^*\left(\frac{\sqrt[3]{9t(1+t+t^2)}}{1+2t}\right)}\leq
\frac{1}{C}\left[\frac{3\pi^2}{4{\sin^2(\pi a)}}\right]
\end{equation}
holds for $a\in[1/3,1/2]$ and $r\in(0,1)$, with equality of (2.22)
or (2.23) if and only if $a=1/3$.

On the other hand, for $a\in(0,1/3]$ ($a\in[1/3,1/2]$ resp.), then
from (1.12) and Theorem 2.1 we have
\begin{align}
\mu_{a}^*(r)\mu_{a}^*\left(\frac{1-r}{1+2r}\right)&=\frac{\pi}{2\sin{(\pi
a)}}\frac{F(a,1-a;1;\frac{{9r(1+r+r^2)}}{(1+2r)^3})}{F(a,1-a;1;\left(\frac{1-r}{1+2r}\right)^3)}\mu_{a}^*(r)\\
&\leq (\geq \mbox{resp.}) \frac{\pi}{2\sin{(\pi
a)}}\frac{(1+2r)F(a,1-a;1;r^3)}{\left(\frac{1+2r}{3}\right)F(a,1-a;1;1-r^3)}\mu_{a}^*(r)\nonumber\\
&=3\mu_{a}^*(r)\mu_{a}^*(r^*)=\frac{3\pi^2}{4{\sin^2(\pi
a)}}.\nonumber
\end{align}
The second equality in (2.24) holds if and only if $a=1/3$. Thus
inequalities (2.15) and (2.16) follows from (1.12) and (2.22)-(2.24)
together with the monotonicity of $f$. $\Box$

\medskip
{\sc Remark.} Theorem 2.7 extends the formulas (1.13) and (1.14) to
the function $\mu_{a}^*(r)$ for $a\in(0,1/2]$.

\bigskip
\section{Proofs of Theorems 1.2 and 1.4}
\setcounter{equation}{0} In this section, we prove our main results
stated in Section 1.

\medskip
{\bf 3.1 Proof of Theorem 1.2.} Consider the function
\begin{equation}
f(r)=\mu_{a}^*(r)+\frac{1}{2}\log\left(\frac{1-r^*}{1+2r^*}\right)
=\mu_{a}^*(r)+\frac{3}{2}\log{r}-\frac{1}{2}\log(1+2r^*)(1+r^*+{r^*}^2)
\end{equation}
for $a\in(0,1/2]$ and $r\in(0,1)$. Let $r_{0}=r^*$,
$r_{1}=\varphi_{3}^*(r^*)=\sqrt[3]{9r^*(1+r^*+{r^*}^2)}/(1+2r^*)$,
$r_{2}=\varphi_{3}^*(r_{1})=\varphi_{9}^*(r^*)$, then
$r^*=\varphi_{1/3}^*(r_{1})=(1-r_{1}^*)/(1+2r_{1}^*)$,
$r_{1}^*=(1-r^*)/(1+2r^*)=\varphi_{1/3}^*(r)$, and
$r=\varphi_{3}^*(r_{1}^*)$ so that
\begin{align}
f(r)=&\mu_{a}^*\left(\frac{\sqrt[3]{9r_{1}^*(1+r_{1}^*+{r_{1}^*}^2)}}{1+2r_{1}^*}\right)
+\frac{1}{2}\log\left(\frac{1-\varphi_{1/3}^*(r_{1})}{1+2\varphi_{1/3}^*(r_{1})}\right)\\
=&\mu_{a}^*\left(\frac{\sqrt[3]{9r_{1}^*(1+r_{1}^*+{r_{1}^*}^2)}}{1+2r_{1}^*}\right)+\frac{1}{2}\log{r_{1}^*}.\nonumber
\end{align}

Let
$g(x)=3\mu_{a}^*\left(\sqrt[3]{9x(1+x+{x}^2)}/(1+2x)\right)-\mu_{a}^*(x)$
for $x\in(0,1)$ and $a\in(0,1/2]$. Then (3.2) can be written as
\begin{equation*}
f(r)-\frac{1}{2}\log{r_{1}^*}+\frac{1}{6}\log\left(\frac{1-r_{1}}{1+2r_{1}}\right)=\frac{1}{3}[g(r_{1}^*)+f(r_{1}^*)],
\end{equation*}
that is
\begin{equation}
f(r)-\frac{1}{6}\log[(1+2r_{1})(1+r_{1}+r_{1}^2)]=\frac{1}{3}[g(r_{1}^*)+f(r_{1}^*)].
\end{equation}

Similarly, putting $r_{2}=\varphi_{3}^*(r_{1})=\varphi_{9}^*(r^*)$,
we get
\begin{equation}
f(r_{1}^*)-\frac{1}{6}\log[(1+2r_{2})(1+r_{2}+r_{2}^2)]=\frac{1}{3}[g(r_{2}^*)+f(r_{2}^*)],
\end{equation}
and hence, by (3.3),
\begin{align}
&f(r)-\frac{1}{6}\log[(1+2r_{1})(1+r_{1}+r_{1}^2)]-\frac{1}{18}\log[(1+2r_{2})(1+r_{2}+r_{2}^2)]\\
=&\frac{1}{3}g(r_{1}^*)+\frac{1}{9}g(r_{2}^*)+\frac{1}{9}f(r_{2}^*).\nonumber
\end{align}

Generally, assuming
\begin{equation}
f(r)-\frac{1}{2}\sum_{k=1}^{n-1}\frac{1}{3^k}\log[(1+2r_{k})(1+r_{k}+r_{k}^2)]=\sum_{k=1}^{n-1}\frac{1}{3^k}g(r_{k}^*)+\frac{1}{3^{n-1}}f(r_{n-1}^*)
\end{equation}
for $n\in \mathbb{N}$ and $n\geq 2$, we let
$r_{n}=\varphi_{3}^*(r_{n-1})=\varphi_{3^n}^*(r^*)$ in (3.3), and
from (3.6) it follows that
\begin{equation}
f(r)-\frac{1}{2}\sum_{k=1}^{n}\frac{1}{3^k}\log[(1+2r_{k})(1+r_{k}+r_{k}^2)]=\sum_{k=1}^{n}\frac{1}{3^k}g(r_{k}^*)+\frac{1}{3^{n}}f(r_{n}^*).
\end{equation}
Hence, by induction, (3.7) holds for all $n\in\mathbb{N}$,
$a\in(0,1/2]$, and $r\in(0,1)$.

Next, we divide the proof into two cases.

{\bf Case 1} $a\in(0,1/3]$. Then from (3.7), Lemma 2.6(2) and
Theorem 2.7(1) we have
\begin{align*}
\frac{1}{3^n}\left(-\frac{3\log{3}}{2}\right)&\leq
f(r)-\frac{1}{2}\sum_{k=1}^{n}\frac{1}{3^k}\log[(1+2r_{k})(1+r_{k}+r_{k}^2)]\nonumber\\
&\leq \sum_{k=1}^{n}\frac{1}{3^k}(R(a)-\log{27})+
\frac{1}{3^n}\frac{R(a)}{2}\nonumber\\
&=\frac{1}{2}(R(a)-\log{27})+\frac{\log{27}}{2\cdot 3^n}.
\end{align*}
Letting $n\rightarrow \infty$, we get
\begin{align}
&\frac{1}{2}\sum_{k=1}^{\infty}\frac{1}{3^k}\log[(1+2r_{k})(1+r_{k}+r_{k}^2)]\leq
f(r)\\
&\leq
\frac{1}{2}\sum_{k=1}^{\infty}\frac{1}{3^k}\log[(1+2r_{k})(1+r_{k}+r_{k}^2)]+\frac{1}{2}(R(a)-\log{27}).\nonumber
\end{align}

The double inequality (1.18) follows from (3.1) and (3.8).

{\bf Case 2} $a\in[1/3,1/2]$. It follows from (3.7), Lemma 2.6 and
Theorem 2.7(1) that
\begin{align*}
\frac{1}{2}(R(a)-\log{27})-\frac{R(a)}{2\cdot 3^n}=&\sum_{k=1}^{n}\frac{1}{3^k}(R(a)-\log{27})+\frac{1}{3^n}\left(-\frac{3\log{3}}{2}\right)\nonumber\\
&\leq
f(r)-\frac{1}{2}\sum_{k=1}^{n}\frac{1}{3^k}\log[(1+2r_{k})(1+r_{k}+r_{k}^2)]\nonumber\\
&\leq \frac{R(a)}{2\cdot3^n}.
\end{align*}
Letting $n\rightarrow \infty$, we get
\begin{align}
&\frac{1}{2}\sum_{k=1}^{\infty}\frac{1}{3^k}\log[(1+2r_{k})(1+r_{k}+r_{k}^2)]+\frac{1}{2}(R(a)-\log{27})\leq
f(r)\\
&\leq
\frac{1}{2}\sum_{k=1}^{\infty}\frac{1}{3^k}\log[(1+2r_{k})(1+r_{k}+r_{k}^2)].\nonumber
\end{align}

The double inequality (1.19) follows from (3.1) and (3.9), and the
remaining results are clear. $\Box$

\medskip
The following corollary follows easily from Theorem 1.2.

\medskip
{\sc Corollary 3.1.} \emph{Let $r\in(0,1)$, then
\begin{equation}
\mu^*(r)\leq \mu_{a}^*(r) \leq \mu^*(r)+\frac{1}{2}(R(a)-\log{27})
\end{equation}
if $a\in(0,1/3]$, and
\begin{equation}
\mu^*(r)+\frac{1}{2}(R(a)-\log{27}) \leq \mu_{a}^*(r) \leq \mu^*(r)
\end{equation}
if $a\in[1/3,1/2]$, each equality in (3.10) or (3.11) is reached if
and only if $a=1/3$.}

\medskip
{\bf 3.2 Proof of Theorem 1.4.} Since the limit function $F(a,b)$
satisfies
\begin{equation*}
F(a_{0},b_{0})=F(a_{1},b_{1})=\cdots,
\end{equation*}
we get
\begin{equation*}
F(a_{0},b_{0})=F(a_{1},b_{1})=F\left(\frac{a_{0}+2b_{0}}{3},\sqrt[3]{\frac{b_{0}(a_{0}^2+a_{0}b_{n}+b_{0}^2)}{3}}\right)
\end{equation*}
or
\begin{align}
F(x,1)=&F\left(\frac{x+2}{3},\sqrt[3]{\frac{x^2+x+1}{3}}\right)\\
=&\sqrt[3]{\frac{x^2+x+1}{3}}F\left(\frac{x+2}{\sqrt[3]{9(x^2+x+1)}},1\right).\nonumber
\end{align}

Let
\begin{equation*}
H(x)=\frac{x^{1/2}(1-x)^{1/3}}{F((1-x)^{1/3},1)},
\end{equation*}
\begin{equation}
t(x)=1-\frac{(\widehat{x}+2)^3}{9({\widehat{x}}^2+\widehat{x}+1)}=\frac{(1-\widehat{x})^3}{9({\widehat{x}}^2+\widehat{x}+1)},
\quad \widehat{x}=(1-x)^{1/3}.
\end{equation}
Then from equations (3.12) and (3.13) we have
\begin{equation}
\frac{d\widehat{x}}{dx}=-\frac{1}{3\widehat{x}^2},\quad
\frac{dt(x)}{dx}=\frac{(1-\widehat{x})^2(\widehat{x}+2)^2}{27{\widehat{x}}^2({\widehat{x}}^2+\widehat{x}+1)^2},
\end{equation}
and
\begin{align}
\frac{H(x)}{H(t(x))}=&\frac{x^{1/2}(1-x)^{1/3}}{\frac{(1-\widehat{x})^{3/2}}{3({\widehat{x}}^2+\widehat{x}+1)^{1/2}}\cdot
\frac{(\widehat{x}+2)}{9^{1/3}({\widehat{x}}^2+\widehat{x}+1)^{1/3}}}\cdot
\sqrt[3]{\frac{3}{{\widehat{x}}^2+\widehat{x}+1}}\\
=&\frac{9\widehat{x}({\widehat{x}}^2+\widehat{x}+1)}{(1-\widehat{x})(\widehat{x}+2)}=\sqrt{\frac{3}{t'(x)}}.\nonumber
\end{align}
Thus the key point of the proof is to show that
\begin{equation}
G(x)=x^{1/2}(1-x)^{1/3}F(\frac{1}{3},\frac{1}{3};1;x)
\end{equation}
also satisfies the function equation (3.15). From this we deduce
that $G(x)=H(x)$. In fact, if we let $J(x)=G(x)/H(x)$, then
$J(x)=J(t(x))$. Note that
\begin{equation*}
t(x)=\frac{(1-\widehat{x})^3}{9({\widehat{x}}^2+\widehat{x}+1)}<\frac{1}{9}x
\end{equation*}
for $x\in(0,1)$, thus $J(x)=J(0^+)=1$ for $x\in(0,1)$.

The hypergeometric differential equation satisfied by $G$ is
\begin{align}
\frac{G''(x)}{G(x)}=\frac{-9x^2+10x-9}{36x^2(1-x)^2}=a(x),
\end{align}
since $F(\frac{1}{3},\frac{1}{3};1;x)$ satisfies the hypergeometric
differential equation
\begin{equation*}
x(1-x)y''+\left(1-\frac{5}{3}x\right)y'-\frac{1}{9}y=0.
\end{equation*}
Now it is a calculation that
\begin{equation}
G^*(x)=\sqrt{\frac{3}{t'(x)}}G(t(x))
\end{equation}
also satisfies (3.17) exactly when
\begin{equation}
a(x)=(t'(x))^2a(t(x))-\frac{1}{2}\frac{t'''(x)}{t'(x)}+\frac{3}{4}\left(\frac{t''(x)}{t'(x)}\right)^2.
\end{equation}

A tedious calculation gives
\begin{equation*}
t''(x)=\frac{2(1-\widehat{x})(\widehat{x}+2)(-{\widehat{x}}^4-2{\widehat{x}}^3+6{\widehat{x}}^2+4\widehat{x}+2)}{81{\widehat{x}}^5({\widehat{x}}^2+\widehat{x}+1)^3},
\end{equation*}
\begin{equation*}
t'''(x)=\frac{2(5{\widehat{x}}^{8}+20{\widehat{x}}^7-34{\widehat{x}}^{6}-136{\widehat{x}}^5-22{\widehat{x}}^4+68{\widehat{x}}^3+104{\widehat{x}}^2+56\widehat{x}+20)}{243{\widehat{x}}^8({\widehat{x}}^2+\widehat{x}+1)^4},
\end{equation*}
\begin{equation*}
a(t(x))=-\frac{81({\widehat{x}}^2+\widehat{x}+1)^2({\widehat{x}}^{6}+4{\widehat{x}}^5+76{\widehat{x}}^4+152{\widehat{x}}^3+248{\widehat{x}}^2+176\widehat{x}+72)}{4(1-\widehat{x})^6({\widehat{x}}^3+6{\widehat{x}}^2+12\widehat{x}+8)}.
\end{equation*}
Putting the three equations above into the right hand of (3.19), by
simplification using {\textsc{Maple 13}}, we get
\begin{align*}
&(t'(x))^2a(t(x))-\frac{1}{2}\frac{t'''(x)}{t'(x)}+\frac{3}{4}\left(\frac{t''(x)}{t'(x)}\right)^2\\
=&-\frac{9{\widehat{x}}^6-8{\widehat{x}}^3+8}{36{\widehat{x}}^6(1-\widehat{x})^2({\widehat{x}}^2+\widehat{x}+1)^2}=-\frac{9(1-x)^2-8(1-x)+8}{36x^2(1-x)^2}=a(x).\nonumber
\end{align*}
Hence both $G^*(x)$ and $G(x)$ satisfy (3.17). Furthermore, since
the roots of the indicial equation of (3.17) are $(1/2,1/2)$ there
is a fundamental logarithmic solution. Since both $G^*$ and $G$ are
asymptotic to $\sqrt{x}$ at $0$, they are in fact equal. Thus (3.18)
shows that $G$ satisfies (3.15). This implies
$F(x,1)=1/F(1/3,1/3;1;1-x^3)$, and equation (1.22) follows from
(3.12) . $\Box$

\end{document}